# Using working patterns as a basis for differentiating part-time employment

**Patrick Letrémy * — Christèle Meilland ** — Marie Cottrell ***

*\* Samos-Matisse, CNRS UMR 8595*
*pley@univ-paris1.fr*
*cottrell@univ-paris1.fr*

*\*\* IRES*
*christele.meilland@ires-fr.org*

ABSTRACT. *Seeking to determine which working patterns have a specific effect on part-time work, in 1998-99 France's INSEE statistical agency carried out a Timetable survey that questioned the homogeneity of this form of employment (again in terms of the working patterns upon which it is based). A neuronal method was used to classify an entire sample of part-time employees according to their weekly working patterns –the end result being that part-time work was shown to be a very heterogeneous form of employment. This was not only reflected by the existence of many different groups of part-time employees, each with highly differentiated individual and professional characteristics, but also (and above all) by the diversity of their weekly working patterns.*

CLASSIFICATION CODE:

## 1. Introduction

In 2001, 16% of France's currently employed working population was working on a part-time basis, versus 7.7% in 1982 (INSEE, 1982 and 2001). Part-time work has experienced unprecedented growth in this country over the past 20 years, notably since 1992 in the wake of a series of State incentives. Clearly this form of employment, often described as being "particular" in nature (as opposed to a norm comprised of a full-time, open-ended employment contract) is still far from representing the numerical majority but it has played an undeniable role in helping women to enter the labour market. 30.4% of currently employed workingwomen (versus 5% of all currently employed persons in 2001) presently exercise their profession on a part-time basis.

But to whom does the expression "part-time worker" actually refer? A part-time working population does not constitute a homogeneous group. After all, employees work on a part-time basis for a variety of reasons or motivations (labour market constraints, family constraints, health constraints, involuntary exit from the labour market, being forced to choose under constraint, lacking any say in one's situation, etc.) – and a whole variety of working time patterns characterise this kind of work. In actual fact, the only factor that the different types of part-time employment have in common (in terms of their working hours) is that people who find themselves in this sort of situation work fewer hours than they would in a full-time job, i.e., their time at work has been shortened. As for their working patterns, these also differ greatly from one another.

INSEE's 1998-99 Timetable survey (Letrémy, Macaire et. al., 2001) was an attempt to determine which working patterns have a specific effect on part-time work; and to question the homogeneity of this form of employment (again, in terms of the working patterns it infers). This explains the survey's use of a neuronal method to classify a population of part-time employees. What this did was to cast doubts as to the homogeneity of part-time work, in terms of the weekly working patterns it involves. The following question was also answered: are part-time jobs predicated upon specific or on standard types of working patterns?

## 2. Data

To analyse weeklong daily working patterns, two data files deriving from the same INSEE/DARES 1998-99 Timetable Survey had to be re-combined. The first (the "individual questionnaire" file) presented interviewees' individual characteristics - the second each individual's professional working patterns over the course of a week.

The goal of the first questionnaire was to identify both individual and professional characteristics (working hours, working patterns, activity levels and profession). An initial study entitled "Working times in particular forms of employment: the specific case of part-time work" (Letrémy & Cottrell, 2001, 2002) covered 14 of the questions that the questionnaire had asked, representing 39 response modalities and 827 part-time workers. The main tool that this study used was the KDISJ algorithm, derived from Kohonen's algorithm (Kohonen 1984, 1993, 1995). Table 1 lists the variables and response modalities that were included in this initial study.

| Heading | Name | Response modalities |
|---|---|---|
| Nature of employment contract | **Contract** | Open-ended contract, fixed-term contract |
| Sex | **Sex** | Man, Woman |
| Age | **Age** | <25, [25, 40[, [40,50[, ≥50 |
| Daily work schedules | **DaySch** | Identical, as-Posted, Variable |
| Number of days worked in a week | **DayWk** | Identical, Variable |
| Night work | **Night** | Usually, sometimes, never |
| Saturday work | **Sat** | Usually, sometimes, never |
| Sunday work | **Sun** | Usually, sometimes, never |
| Wednesday work | **Wed** | Usually, sometimes, never |
| Ability to go on leave | **Leave** | Yes no problem, yes under conditions, no |
| Who determines employee's schedule | **Det** | Company, a la carte, employee, other |
| Involuntary nature of part-time status | **Volunt** | Yes (involuntary), no (voluntary) |
| Awareness of next day's schedule | **Next** | Yes, no |
| Possibility of carrying over credit hours | **Carry** | No point, yes, no |

**Table 1.** *Variables that were used in the individual survey*

In the second questionnaire (the weekly report), the goal was to ascertain each individual's daily and weekly working patterns on a quarter-hourly basis. Every day interviewees would fill in a sheet stating for each quarter-hour whether they actually worked, i.e., respondents marked (1) if they had worked or (0) if they had not for 7 days in a row, accounting for a total of 4 x 24 x 7 = 672 quarter-hours. The sum total of these responses constituted the "weekly report". Each individual would then be attributed a weekly working profile, constituted on the basis of binary values.

These working patterns were observed over a continuous week. Interviewees provided responses about their schedules and working times for 7 days in a row before mailing the completed questionnaire back to INSEE. This procedure is worthy of mention since it explains why certain weekly reports went missing and why others were not entirely accurate.

Due to these incomplete or non-existent weekly reports, the study ended up covering fewer people[1] once the two disjointed data files had been recombined using statistical and/or IT techniques. Its working patterns analysis only dealt with an overall population of part-time employees, irrespective of whether said individuals were on an open-ended or a fixed-term employment contract. A total of 566 employees were studied, broken down into 473 open-ended employment contract holders vs. 93 fixed-term contract holders. The population was mixed but not divided equally, with 505 women vs. 61 men.

All in all, two sorts of data were used: the first type in the shape of tables and test values derived from data contained within the individual questionnaire; and the second in the shape of professional occupation profiles derived from the data contained within the weekly reports.

In most cases, the two data files featured similar sorts of outcomes, with the results of the individual file confirming analysis derived from the weekly report (or vice versa). Occasionally however the responses did diverge. These divergences and contrasts should be viewed as reflecting the nature of these questionnaires; and how difficult it is to weld distinctive forms of data files.

For example, working times were described in the individual Timetable survey in at least three different ways. One direct estimate came out of the question "How many hours do you actually work in a normal week?", with interviewees having to indicate a maximum and a minimum. Another came from the question "Theoretically how much time are your supposed to work every week (in hours and minutes)?" Total working hours could be calculated in the weekly report on the basis of the number of quarter-hours of work [2] done in a supposedly normal week.

Similarly, information about working nights, Saturdays, Sundays or Wednesdays (Note Trans. - many French children only go to school for half a day on Wednesdays, and this has a effect on some people's work schedules) can be found in the individual survey (in answers to the Night, Sat, Sun and Wed questions) but clearly also in the corresponding days and hours that were specified in the weekly report. There were also a few incoherent results that we will attempt to explain.

---

[1] As opposed to the initial population (Letrémy & Cottrell, 2001).
[2] Professional work brought home is included but time for meals and commuting time is not.

## 3. Weekly report rankings

We will not be reminding you here of how Kohonen's algorithm can be defined or applied to data analysis (see for example Cottrell et. al.; 1998, Cottrell & Rousset, 1997; Kaski, 1997).

A 10-unit Kohonen string (a uni-dimensional map) is used [3] to classify 566 weekly profiles formed from binary vectors in 672 dimensions. Each class is then represented by (summarised in) a 672-dimensional code vector. The x-axis is the time (quarter-hour during the day), starting Monday 0:00 and continuing until Sunday midnight. The y-axis is a number included between 0 and 1, obtained for each quarter-hour of each day of the week. This can be interpreted as the proportion of individuals in the class who are considered to be actively working at that moment in time.

Figure 1 highlights these 10 code vectors, using a vertical line to separate each day. At first glance, we can see right away that the 10 classes are very clear and distinct, and that the code vectors are perfectly ordered from top ("normal" working conditions) to bottom ("non-standard" working conditions).

We then recombine these 10 classes by classifying the 10 code vectors hierarchically. We do this until we get to the point where 85.6% of the variance is accounted for by five superclasses. This corroborates the Kohonen map's organisational quality, insofar as only those classes that are consecutive will be recombined, and this on a two-by-two basis. These five superclasses are marked A, B, C, D and E. Each contains two code vectors: class A (1 and 2), class B (3 and 4), class C (5 and 6), class D (7 and 8) and class E (9 and 10).

The 5 superclasses' sample sizes are fairly evenly balanced:

|             | A   | B   | C   | D   | E   |
|-------------|-----|-----|-----|-----|-----|
| Sample size | 141 | 100 | 108 | 110 | 107 |

---

[3] All of the software used here was written by Patrick Letrémy in SAS and is available on the SAMOS website: http://samos.univ-paris1.fr)

Figure 1 represents the 10 classes with their code vectors. The column on the left indicates their recombination into superclasses.

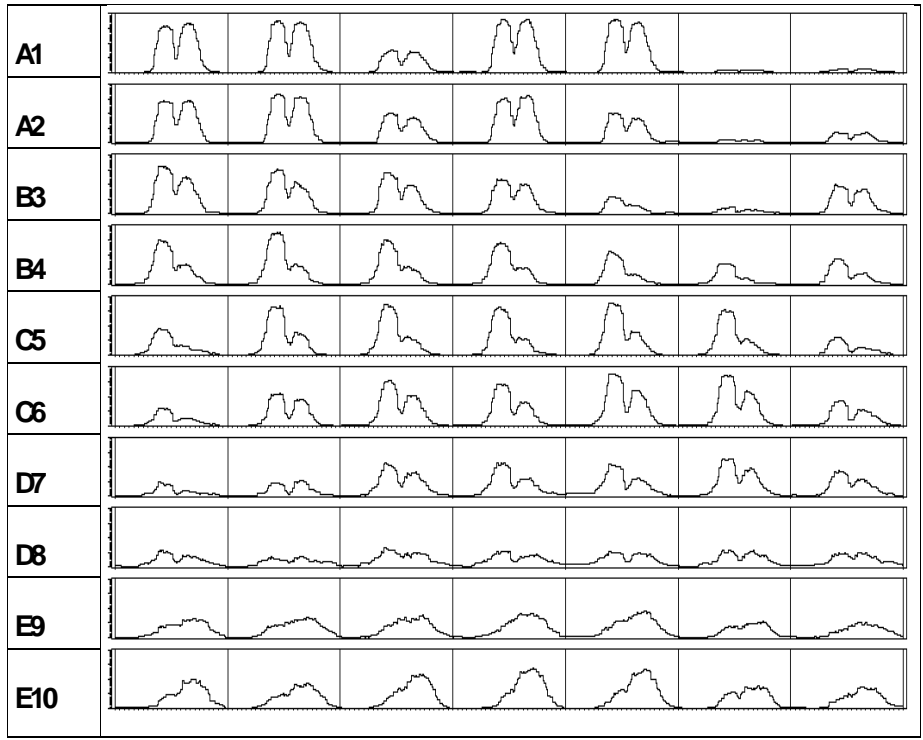

**Figure 1.** *Profiles of the 10 weekly report classes (from Monday through Sunday)*

These profiles can be organised "correctly" using a "Multi Dimensional Scaling" (MDS) technique. This gives us the uni-dimensional structure of a Kohonen classification over a string of 10 units.

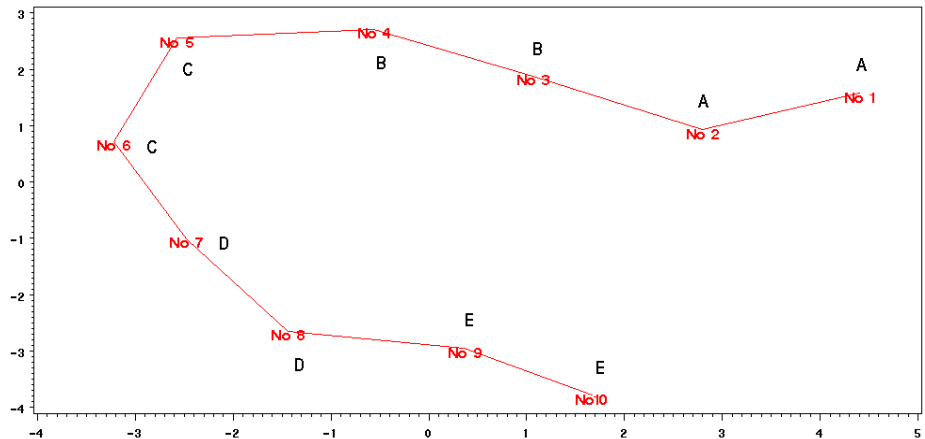

**Figure 2** *Representation of the 10 code vectors or typical profiles, after using an MDS technique. Note the perfect organisation, i.e., there is no crossing anywhere within the string*

Since we know how to identify the individuals who belong to each of the five superclasses, we can cross the classification we have obtained in this manner with the 14 questions featured in the individual questionnaire. We start by using chi-two independence tests to highlight which questions are discriminatory (in terms of our partition into five classes). We then remove three of these questions, those where tests were not significant (where the A, B, C, D and E class modalities and the question modalities were independent). In other words, we do not take gender, number of days per week and awareness of the following day's schedule into consideration.

Amongst the 11 questions that we do consider, two levels of description can be distinguished:

1) An initial level supplementing the information provided by the weekly report profiles that typify classes A, B, C, D and E. We only consider the seven questions (Contract, Age, Daily work schedules (DaySch), Ability to go on leave (Leave), Who determines schedules (Det), Involuntary nature of part-time status (Volunt) and Possibility of carrying over credit hours (Carry) that pertain to age and working conditions.

2) The second level checks the coherency of the information provided in the weekly report charts or in the questionnaire. This includes four questions relative to working nights (Night), Saturdays (Sat), Sundays (Sun) and Wednesdays (Wed). As aforementioned, there can be some divergence between an individual's response to the questionnaire and what actually transpired over the seven days that s/he was filling in the weekly report sheet. To assess the significance of this problem, we counted the number of individuals actually at work at specific times of the day, i.e., 10AM, 4PM and 9PM.

This partial crossing of variables focuses more specifically on percentages of response modalities which had a test value that was strictly greater than 1 (meaning where the response rate was greater for this one class than it was for the whole of the survey). This has led to a description that is both based on the five superclasses A, B, C, D and E and also separated into two levels.

| QUESTION | Categories | A | B | C | D | E | Total |
|---|---|---|---|---|---|---|---|
| **Contract** | Open-ended | **89** | **87** | *81* | *83* | *77* | 84 |
|  | Fixed-term | *11* | *13* | **19** | **17** | **23** | 16 |
| *Age* | < 25 | *4* | *4* | *2* | **9** | **20** | 8 |
|  | [25, 40 [ | **40** | *33* | **42** | **40** | **39** | 39 |
|  | [40, 50 [ | **33** | **35** | **31** | *29* | *26* | 31 |
|  | >= 50 | **23** | **28** | **25** | *22* | *15* | 22 |
| *DaySch* | Identical | **52** | **61** | **59** | *47* | *36* | 51 |
|  | as-Posted | *1* | **5** | *2* | **5** | **7** | 4 |
|  | Variable | **46** | *34* | *39* | **47** | **57** | 45 |
| *Leave* | Yes | **77** | **76** | *69* | *75* | **76** | 75 |
|  | Yes but | **15** | **16** | **18** | *13* | *7* | 14 |
|  | No | *9* | *8* | **13** | **13** | **18** | 12 |
| *Det* | Company | *52* | **64** | *58* | **65** | **75** | 62 |
|  | A la carte | *9* | *7* | **17** | **12** | *6* | 10 |
|  | Employee | **36** | **21** | *19* | *14* | *7* | 20 |
|  | Others | *3* | **8** | *6* | **9** | **12** | 7 |
| *Volunt* | Involuntary | *35* | *51* | *47* | **62** | **65** | 51 |
|  | Voluntary | **65** | *49* | **53** | *38* | *35* | 49 |
| *Carry* | NA | *50* | **58** | *53* | **59** | *53* | 54 |
|  | Yes | **27** | **27** | **28** | *26* | *26* | 27 |
|  | No | **23** | *15* | **19** | *15* | **21** | 19 |

**Table 2:** *First level. Percentages associated with a test value strictly greater than 1 are in bold font*

This has led to the following deduction(s):

---

**A**: Open-ended contract; part-time work is voluntary; leave granted right away or with just a few conditions; schedule chosen a la carte; age: 40% in [25, 40[, 33% in [40, 50[ and 23% over 50.

**B**: Open-ended contract; part-time work is involuntary; 61% of daily schedules are identical vs. 5% as-posted; leave provisions similar to A; the company determines the schedule although 21% are chosen a la carte; for 58% there is no point in carrying over credit hours but for 27% this is a possibility; age: 35% in [40, 50[ and 28% over 50 (this is an older class than A).

**C**: Fixed-term contract; part-time work is voluntary; daily schedules are identical; for 18% leave is possible with a few conditions but for 13% this is impossible; for 17% schedules can be arranged; age: 42% in [25, 40[, 31% in [40, 50[ and 25% over 50 (comparable to A).

**D**: Fixed-term contract; part-time work is involuntary; for 47% daily schedules are variable but for 5% they are as-posted; leave is impossible; the company determines the schedule for 65% but 12% have a possible choice; no point in carrying over credit hours ; age: 9% under 25 and 40% in [25, 40[.

**E**: Fixed-term contract; part-time work is involuntary; for 57% daily schedule is variable but for 7% it is as-posted; leave is possible for 76% and impossible for 18%; the company determines the schedule; no carryover of credit hours; age: 20% under 25 and 39% in [25, 40[.

---

This typology can be rounded out by a second-level description monitoring the coherency of people's experience as related in the weekly report (headcount at 10AM, 4PM and 9 PM Saturdays, Sundays and Wednesdays) and the answers provided to the four questions relative to working nights, Saturdays, Sundays and Wednesdays.

| In% | Never | Sometimes | Usually | Number of persons working at 10h, 16h and 21h | | | | | | | | |
|---|---|---|---|---|---|---|---|---|---|---|---|---|
| **A** | | | | **A** | | *%* | | *%* | | | | *%* |
| Night | **91** | **7** | *1* | Total | 141 | | | | | | | |
| Sat | **72** | **23** | *6* | Sat_10h | 8 | *6* | Sat_16h | 4 | *3* | Sat_21h | 0 | *0* |
| Sun | **83** | *15* | *2* | Sun_10h | 8 | *6* | Sun_16h | 6 | *4* | Sun_21h | 1 | *1* |
| Wed | **28** | **21** | *51* | Wed_10h | 57 | *40* | Wed_16h | 55 | *39* | Wed_21h | 1 | *1* |
| **B** | | | | **B** | | | | | | | | |
| Night | **94** | *4* | *2* | Total | 100 | | | | | | | |
| Sat | **60** | **21** | *19* | Sat_10h | 12 | *12* | Sat_16h | 8 | *8* | Sat_21h | 3 | *3* |
| Sun | **77** | **20** | *3* | Sun_10h | 60 | *60* | Sun_16h | 38 | *38* | Sun_21h | 2 | *2* |
| Wed | **25** | *7* | **68** | Wed_10h | 79 | *79* | Wed_16h | 35 | *35* | Wed_21h | 2 | *2* |
| **C** | | | | **C** | | | | | | | | |
| Night | **94** | *6* | *0* | Total | 108 | | | | | | | |
| Sat | *41* | *19* | **40** | Sat_10h | 100 | *93* | Sat_16h | 49 | *45* | Sat_21h | 1 | *1* |
| Sun | **78** | *13* | **9** | Sun_10h | 32 | *30* | Sun_16h | 21 | *19* | Sun_21h | 2 | *2* |
| Wed | *20* | *12* | **68** | Wed_10h | 88 | *81* | Wed_16h | 49 | *45* | Wed_21h | 2 | *2* |
| **D** | | | | **D** | | | | | | | | |
| *Night* | *79* | *10* | *11* | *Total* | *110* | | | | | | | |
| Sat | *35* | *18* | **46** | Sat_10h | 43 | *39* | Sat_16h | 34 | *31* | Sat_21h | 4 | *4* |
| Sun | *67* | **19** | **14** | Sun_10h | 40 | *36* | Sun_16h | 25 | *23* | Sun_21h | 8 | *7* |
| Wed | **25** | *13* | **63** | Wed_10h | 45 | *41* | Wed_16h | 25 | *23* | Wed_21h | 8 | *7* |
| **E** | | | | **E** | | | | | | | | |
| Night | **93** | *5* | *2* | Total | 107 | | | | | | | |
| Sat | *30* | **21** | **49** | Sat_10h | 19 | *18* | Sat_16h | 32 | *30* | Sat_21h | 9 | *8* |
| Sun | *76* | *16* | **8** | Sun_10h | 10 | *9* | Sun_16h | 29 | *27* | Sun_21h | 10 | *9* |
| Wed | *16* | **20** | **64** | Wed_10h | 20 | *19* | Wed_16h | 51 | *48* | Wed_21h | 13 | *12* |

**Table 3:** *Second level. The values associated with test values greater than 1 are in bold font. The right-hand side of the table indicates headcounts at 10h (10 AM), 16h (4 PM), 21h (9 PM), Saturdays, Sundays and Wednesday*

The class descriptions can be rounded out as follows:

> **A**: Work neither nights, Saturdays or Sundays, and have a lower activity level on Wednesdays. The weekly report and the individual questionnaire converge.
>
> **B**: Do not work nights; very limited activity levels Saturdays (slightly less than in the questionnaire); mostly active in the morning the other days of the week. With respect to the Sunday question, a clear divergence between the weekly report (60% are at work at 10AM) and the questionnaire (77% state they never work Sundays).
>
> **C**: Do not work nights; mostly work Wednesday and Saturday morning with lower level of activity Sundays (and Mondays). With respect to the Sunday question, a slight divergence between the weekly report (30% are at work at 10AM) and the questionnaire (9% state they usually work Sundays).
>
> **D**: A little night work but less than in the questionnaire; reduced activity level Saturdays and Sundays; mainly work Wednesday mornings (much lower activity levels Mondays and Tuesdays). Little divergence between the weekly report and the questionnaire.
>
> **E**: A little night work (but more than in the questionnaire); mostly work Monday or Friday afternoons with lower activity levels Saturdays and Sundays; slight divergence for Wednesday between weekly report (48% are at work at 4PM) and the questionnaire (64% state they usually works Wednesdays).

## 4. Conclusion

The whole typology can be summarised by representing the average activity levels of the individuals found in each of the five classes. These curves supplement the partial headcounts carried out Wednesdays, Saturdays and Sundays at 10 AM, 4 PM and 9 PM. They mesh perfectly with the typical weekly profiles as determined via the Kohonen classification.

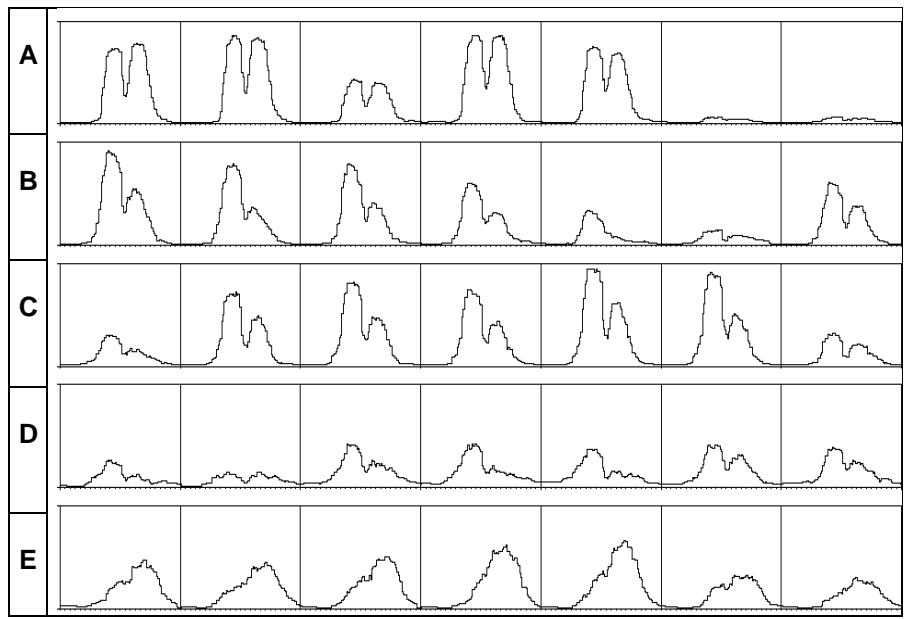

*Figure 3: Average activity level profiles of individuals from each of the five superclasses to have been determined*

The five classes are classified in descending order, both in terms of age and also the quality of working conditions. Moving from A to E they go from a category of open-ended contracts for people over the age of 25 who had volunteered for part-time work and do not work nights or Saturdays and Sundays - to a category comprised of fixed-term contract holders for whom part-time work is involuntary; who are younger (20% are below the age of 25); whose schedules are variable and sometimes only as-posted; and who work Saturdays, Sundays and occasionally nights.

All in all, not only does a summative analysis of the weekly reports (when coupled with a statistical analysis of the individual and professional characteristics of the employees questioned) highlight how very heterogeneous part-time work really is in terms of the working patterns it entails, but also (and even more importantly) it stresses the parallelism that can be traced between the constant rise in situational insecurity between the A and E superclasses (including as-posted or variable schedules, frequent night work and as often as not weekend work – plus less secure conditions of employment, i.e., fixed-term contracts) and one typology of female populations (ranging from the oldest to the youngest). The scale of insecurity we are facing therefore appears to be based on the fact that "voluntary" part-time work often involves open-ended contracts characterised by flexible schedules, or "a la carte" solutions employees have chosen themselves and which often entail regular working patterns and a great deal of freedom to go on leave – these being configurations that are often associated with a female population of medium or advanced age. Inversely (and unsurprisingly enough) "involuntary" work is associated with fixed-term contracts where as-posted basis or variable schedules are often determined by the company; where employees have limited freedom to go on leave; and where people frequently work nights or weekends. The population of "involuntary part-timers" is a young one (9 and 20% of the individuals in the D and E type weeks are below 25). We can presume that this latter situation is a case of people having opted for part-time work "due to a lack of anything better" (Maruani,1996).

Part-time work is therefore pluralistic in nature, in terms of its working patterns; working conditions (night or weekend work); and the population it affects (all women, mostly young). Moreover, this pluralistic trait translates a multitude of motives (market constraints, family constraints, personal preferences, partial retirement, etc.) that induce employees to resort to this form of employment.

Our thanks to Colin Marchika and Alain Chenu (CREST, Laboratory of Quantitative Sociology) for having given us access to the SAS programme so we could read the information in the weekly report.

## 5. Bibliography/References